\documentclass[12pt,a4paper, twoside]{amsart}

\usepackage[margin=2.85cm]{geometry}

\usepackage[T1]{fontenc}
\usepackage[latin1]{inputenc}
\usepackage{amsmath}
\usepackage{amsthm}
\usepackage{amssymb}
\usepackage{hyperref}
\usepackage[UKenglish]{babel}
\usepackage{enumerate}
\usepackage[dvips]{graphicx}

\theoremstyle{definition}

\theoremstyle{remark}
\newtheorem*{rem}{Remark}
\theoremstyle{plain}
\newtheorem{thm}{Theorem}
\newtheorem{lem}[thm]{Lemma}

\newtheorem{cor}[thm]{Corollary}

\newtheorem*{pro}{Problem 1}

\newcommand{\be}{\beta}
\newcommand{\al}{\alpha}
\newcommand{\la}{\lambda}
\newcommand{\F}{\mcal F}

\newcommand{\mb}[1]{\mathbb{#1}}
\newcommand{\mbf}[1]{\mathbf{#1}}
\newcommand{\mcal}[1]{\mathcal{#1}}

\newcommand{\W}{W(\la,  \psi)}
\newcommand{\Wc}{W(\la_0,  \psi)}
\newcommand{\Wz}{W(0,  \psi)}
\newcommand{\w}{W(\mcal C, \la, \psi)}
\newcommand{\ws}{S(\mcal C, \la, \psi, N)}
\newcommand{\wrp}{S^{(I)}(\mcal C, \la, \psi, N)}
\newcommand{\wds}{S^{(II)}(\mcal C, \la, \psi, N)}
\newcommand{\wdsa}{S^{(IIa)}(\mcal C, \la, \psi, N)}
\newcommand{\wdsb}{S^{(IIb)}(\mcal C, \la, \psi, N)}
\newcommand{\wdsba}{S^{(IIb1)}(\mcal C, \la, \psi, N)}
\newcommand{\wdsbb}{S^{(IIb2)}(\mcal C, \la, \psi, N)}
\newcommand{\wdsbbg}{S^{(IIb2)}_G(\mcal C, \la, \psi, N)}
\newcommand{\wdsbbb}{S^{(IIb2)}_B(\mcal C, \la, \psi, N)}
\newcommand{\wdsbbbi}{S^{(IIb2)}_{B(i)}(\mcal C, \la, \psi, N)}
\newcommand{\wdsbbbii}{S^{(IIb2)}_{B(ii)}(\mcal C, \la, \psi, N)}

\renewcommand{\le}{\leqslant}
\renewcommand{\ge}{\geqslant}

\newcommand{\RR}{\mathbb{R}}

\newcommand{\ZZ}{\mathbb{Z}}

\newcommand{\NN}{\mathbb{N}}

\newcommand{\CCC}{\mathcal{C}}

\newcommand{\MMM}{\mathcal{M}}
\newcommand{\HHH}{\mathcal{H}}

\begin{document}

\title[A Jarn\'ik type theorem]{An Inhomogeneous Jarn\'ik type theorem for planar curves}

\author{Dzmitry Badziahin$^\dag$}
\author{Stephen Harrap$^\dag$}
\thanks{$^\dag$ Research supported by EPSRC grant number $EP/L005204/1$l}

\address{D. Badziahin \& S. Harrap, Durham University, Department of Mathematical Sciences, Science Laboratories, South Rd, Durham,
DH1 3LE, United Kingdom}
\email{dzmitry.badziahin@durham.ac.uk, \quad s.g.harrap@durham.ac.uk}
\author{Mumtaz Hussain$^*$ }
\thanks{ $^*$ Research supported by the Australian research council}

\address{M. Hussain, School of Mathematical and Physical Sciences, The University of Newcastle,
Callaghan, NSW 2308, Australia.}
\email{mumtaz.hussain@newcastle.edu.au}

\begin{abstract}
In metric Diophantine approximation there are two main types of approximations: simultaneous and dual for both homogeneous and inhomogeneous settings. The well known measure-theoretic theorems of Khintchine and Jarn\'ik are fundamental in these settings. Recently, there has been substantial progress towards establishing a metric theory of Diophantine approximations on manifolds. In particular, both the Khintchine and Jarn\'ik type results have been established for planar curves except for only one case. In this paper, we prove an inhomogeneous Jarn\'ik type theorem for convergence on planar curves and in so doing complete the metric theory for both the homogeneous and inhomogeneous settings for approximation on planar curves.

\end{abstract}

\maketitle

\section{Introduction}

Classical metric Diophantine approximation deals quantitatively with the density of the rational numbers within the real numbers.
The higher dimensional theory of Diophantine approximation for systems of $m$ linear forms in $n$ variables  combines two different types of classical Diophantine approximation: simultaneous and dual.
Simultaneous Diophantine approximation comprises the component-wise approximation of points $\mbf y=(y_1,\dots, y_m)^T\in \mb R^m$ by $m$-tuples of rational numbers $\{\mbf p/q:(\mbf p, q)\in \mb Z^n\times \mb N\}$, whereas dual Diophantine consists of the approximation of points $\mbf x=(x_1,\dots, x_n)\in \mb R^n$ by `rational' hyperplanes of the form $a_1x_1+\cdots+a_nx_n=a_0$, where $(a_0, \mbf a)= (a_0, a_1, \ldots a_n)\in\mb Z\times \mb Z^n\setminus\{\mbf 0\}$. Over the last few years rich and intricate metric theories have been established relating to both of these types of approximation. Broadly speaking, both theories have followed similar paths of development with advances in the dual theory often following quickly from corresponding breakthroughs in the simultaneous theory, although the methods required are on occasion quite different. In this paper we are concerned with a problem in dual approximation whose development was until recently very far behind the simultaneous analogue. We are able to bring the state of the dual theory in line with the simultaneous theory in this case.

\subsection{Inhomogeneous dual approximation}
Throughout, $\psi:\mb N\to \mb R^+$ will denote a real positive decreasing arithmetic function such that $\psi(n)\to 0$ as $n\to \infty$. We will refer to $\psi$ as an \emph{approximating function} for reasons that will become apparent. Given an approximating function $\psi$ and a real function $\la:\mb R^n\to \mb R$ (which will be referred to as the \textit{inhomogeneous function}), define the set
\begin{equation*}
  \W:=\left\{\mbf x=(x_1,\dots,x_n)\in\mb R^n:\begin{array}{l}
  |a_0+a_1x_1+\cdots+a_nx_n+\la(\mbf x)|<\psi(|\mbf a|) \\[1ex]
  \text{for} \ \ i.m. \ (a_0, \mbf a)\in\mb Z\times \mb Z^n\setminus\{\mbf 0\}
                           \end{array}
\right\},
\end{equation*}
where  $`i.m.$' stands for `infinitely many' and $|\mbf a|=\max\{|a_1|,\dots, |a_n|\}$ is the standard
supremum norm. Any vector $\mbf x \in \RR^n$ will be called {\em $\psi-$approximable (with respect to $\lambda$)} if it lies in the set $\W$. Additionally, in the special case that $\psi(r)=r^{-\tau}$ for
some $\tau>0$ we say that $\mbf x$ is {\em $\tau-$approximable (with respect to $\lambda$)} and denote $\W$ by $W(\la, \tau)$.

When the inhomogeneous function satisfies $\la\equiv 0$ the set $W(\psi):=W(0, \psi)$ represents the classical homogeneous set of $\psi$-approximable vectors within the theory of dual Diophantine approximation. On the other hand, when $\la$ is constant the set $\W$ represents what would usually be referred to as the classical  inhomogeneous set of $\psi-$approximable vectors (in the context of dual approximation). In both of these cases the metric theory associated with the set $\W$ is well understood. In particular, the following is a modern version of a combined statement of three classical theorems of Khintchine (\cite{K24}, 1924),  Groshev (\cite{Groshev}, 1926) and Jarn\'ik (\cite{Jar31}, 1931) and is due to many contributions by many authors at various stages (for example, see \cite{BBDV, BDV_mtl, BVclassical}). It  provides an elegant but simple criterion for the `size' of $\W$ expressed in terms of $s$-dimensional Hausdorff measure $\mcal H^s$. For the definitions of Hausdorff measure and dimension we refer the reader to \S\ref{hm}.
\begin{thm}[Khintchine-Groshev-Jarn\'ik]\label{thm2}
    Let $\psi$ be an approximating function, let $\la(\mbf x)=\lambda_0\in \mb R$ be constant and let $s\in (n-1, n]$. Then,
    $$\mcal H^s(\Wc)=\begin {cases}
    0 \ & {\rm if } \quad \sum\limits_{q=1}^{\infty}\psi(q)^{s-n+1}q^{2n-s-1}<\infty. \\[3ex]
    \mcal H^s(\RR^n) \ & {\rm if } \quad \sum\limits_{q=1}^{\infty}\psi(q)^{s-n+1}q^{2n-s-1}=\infty.
    \end {cases}$$
\end{thm}
Note that when $s = n$, the measure $\mcal H^n$ is equivalent to $n$-dimensional Lebesgue measure.  The convergent case of Theorem \ref{thm2} is actually an easy consequence of the Borel-Cantelli lemma from probabiliy theory; the main substance of the theorem lies in the much more difficult divergent part.
Before we proceed, it is worth mentioning that the
  analogue of Theorem \ref{thm2} has been obtained in the setting of simultaneous Diophantine approximation. Here, instead of $\W$ one may consider the set
\begin{equation}\label{vb2}
  \mcal S(\la, \psi):=\left\{\mbf y\in\mb R^m:\begin{array}{l}
  \max\limits_{1\le i\le m}|qy_i -p_i + \la(\mbf y)|<\psi(q) \\[2ex]
  \text{for} \   i.m. \ (q, p_1,\dots,p_n)\in\mb N\times \mb Z^n
                           \end{array}
\right\}
\end{equation}
for some inhomogeneous function $\lambda: \RR^m \to \RR$ and an approximating function $\psi$. More generally, one can combine the two viewpoints (dual and simultaneous) and in the case that $\lambda$ is constant establish a result akin to Theorem \ref{thm2} for the resulting system of $m$ linear forms in $n$ variables-- see \cite{BBDV} for a detailed account. Due to the path of development, results concerning the Lebesgue measure of $\psi$-well approximable points are often referred to as being of \textit{Khintchine type} for the setting of simultaneous Diophantine approximation and of \textit{Groshev type} for the setting of dual Diophantine approximation. On the other hand, results concerning the Hausdorff measure of such sets are often referred to as being of \textit{Jarn\'ik type}.

In this paper we consider a problem in the setting of `functional' inhomogeneous dual approximation; that is, the setup where the inhomogeneous function $\lambda$ is not constant and so depends on each point~$\mbf x \in \RR^n$ being approximated. Measure-theoretic considerations corresponding to this difficult variant of classical dual Diophantine approximation do not seem to have received any explicit treatment in the literature, even for the case~`$s=n$' corresponding to Lebesgue measure. However, by appealing to modern techniques we are able to provide a precise criterion for the Hausdorff measure of the part of the set~$\W$ lying on any sufficiently well behaved planar curve. 
%
%
%

\subsection{Diophantine approximation on manifolds}

The problem of estimating the size of the set of $\psi$--approximable points is more intricate if one restricts $\mbf x\in\mb R^n$ to lie on some $m$--dimensional submanifold $\mcal M\subseteq \mb R^n$.
 Such a restriction implies that the points of interest $\mbf x$ must be functionally related. Consequently, problems of this type were classically referred to as `approximation of dependent quantities'.  Naturally, one must appeal to the induced $m$-dimensional Lebesgue measure $\Lambda$ on the manifold in question, for otherwise the relevant measure-theoretic statements would be trivial - when $m<n$ the $n$-dimensional Lebesgue measure of $\mcal M\cap \W$ is always zero, irrespective of the approximating function $\psi$.

 In what follows, as usual, $C^{(n)}(I)$ will
 denote the set of $n$-times continuously differentiable functions defined on some interval $I$
 of $\mb R$.

 \subsubsection{The classical results}

The exploration of (homogeneous) Diophantine approximation on manifolds dates back to a profound conjecture of  Mahler \cite{Mah32} in 1932, which can be rephrased in terms of the `extremality' of the Veronese curve
 $\mcal V_n:=\{(x,x^2, \ldots ,x^n):x\in\mb R\}$. A manifold $\mcal M \subset \mb R^n$ is said to be  {\em extremal} if $\Lambda(\Wz\cap \mcal M)=0$ for any  $\tau>n$, and Mahler's conjecture was precisely the statement that~$\mcal V_n$ is extremal.
 One should observe that from a measure-theoretic perspective the property of being extremal is actually a weaker property than the convergent part of Theorem~\ref{thm2}, although of course it embodies the added difficulty that the points $\mbf x$ must be dependent upon one another.

 Mahler's conjecture was solved completely by Sprind\v zuk \cite{Spr_Mah} in 1965, although the
 special cases $n=2$ and $n=3$ had been settled earlier. Schmidt \cite{Sch_cur} extended Mahler's problem to the case of all sufficiently curved $ C^{(3)}$ planar curves, leading to a reasonably
 general theory of Diophantine approximation on manifolds. To be precise, Schmidt required that for each local parametrisation $\mathcal C=\mathcal C_f:=\{(f_1(x), f_2(x)): \; x\in I\}$, where $I \subset \RR$ is some closed interval, the curve $\mathcal C_f$ is \textit{non-degenerate}; that is, the set of points $x \in I$ at which $f_1^\prime(x)f_2^{\prime\prime}(x) - f_1^{\prime\prime}(x)f_2^{\prime}(x) =0$ is a set of Lebesgue measure zero. Following on from this, Sprind\v zuk conjectured that any analytic manifold satisfying a similar `non-degeneracy' condition should also be extremal. The heuristic meaning of non-degeneracy in this context is that each manifold is smooth and `curved enough' so as to deviate from any given hyperplane;  see \cite{Ber99,KM} for a precise formulation. Even though particular cases of Sprind\v zuk's conjecture were known, it was not until 1998 when it was solved in full generality by Kleinbock \& Margulis \cite{KM}. Building on this breakthrough, progress has subsequently been dramatic. In particular, the following Lebesgue measure analogue of Theorem \ref{thm2} for Diophantine approximation on manifolds was proved iteratively in the papers \cite{Ber99, BBKM02, BKM01}, whilst \cite{Ber99} gave an alternative proof of Sprind\v zuk's conjecture.
 \begin{thm}[Beresnevich-Bernik-Kleinbock-Margulis]\label{thm3}
Let $\MMM \subset \RR^n$ be a non-degenerate manifold of dimension $m$ and let $\psi$ be an approximating function. 
Then,
$$\Lambda(\Wz\cap \MMM)=\begin {cases}
 0 \ & {\rm if } \quad \sum\limits_{q=1}^{\infty}\psi(q)q^{n-1} \, < \, \infty. \\[2ex]
 \Lambda(\MMM) \ & {\rm if } \quad \sum\limits_{q=1}^{\infty}\psi(q) q^{n-1} \, = \, \infty.
\end {cases}$$
%
%
\end{thm}

Unlike Theorem~\ref{thm2}, the convergence case of Theorem~\ref{thm3} (which was independently proved in \cite{Ber99} and \cite{BKM01}) required very delicate covering and counting arguments
to reduce it to a situation where the Borel-Cantelli lemma was applicable; this was a highly non-trivial task. The divergence part was first established for the Veronese curves in \cite{B} and later for arbitrary non-degenerate manifolds in \cite{BBKM02}. The statement of Theorem~\ref{thm3} was first conjectured by Baker in 1975.

Obtaining a full Hausdorff measure analogue of Theorem~\ref{thm2} for Diophantine approximation on manifolds represents a deep open problem and is often referred to as the \textit{Generalized Baker-Schmidt Problem for Hausdorff measure}. However, whilst establishing a general framework for attacking problems of this type,  Beresnevich, Dickinson \& Velani were remarkably able to verify the conjecture in the case of divergence (see Theorem~18 of \cite{BDV_mtl}) leading to the following unified statement.

 \begin{thm}[Beresnevich-Dickinson-Velani]\label{thm3a}
    Let $\MMM \subset \RR^n$ be a non-degenerate manifold of dimension $m$ and let $\psi$ be an approximating function. Then, for any $s \in (m-1,m]$ we have
 \begin{equation*}
  \mcal H^s(\Wz\cap\MMM) \: = \: \mcal H^s(\MMM) \quad {\rm if } \quad \sum\limits_{q=1}^{\infty}\psi(q)^{s-m+1}q^{n+m-s-1}=\infty.
 \end{equation*}
\end{thm}
A convergence counterpart of Theorem~\ref{thm4} currently remains out of reach. Indeed, it is only very recently that any progress has been made. In 2014 it was finally shown by Huang~\cite{Huang} that the convergence part of the Generalized Baker-Schmidt Problem for Hausdorff measure holds for all $\CCC^{(2)}$ non-degenerate  planar curves. The only previously known result of this type was for the special case of parabolas \cite{Hus_parabola}.


\subsubsection{`Functional' inhomogeneous Diophantine approximation on manifolds}

The most difficult variant of dual Diophantine approximation on manifolds concerns functional inhomogeneous approximation. Until a short time ago the theory in this setting has been almost non-existent, even in the case where the inhomogeneous function $\lambda$ is constant and non-zero; i.e., even in the setting of classical inhomogeneous approximation. The concept of inhomogeneous extremality  (for $\lambda$ taking a non-zero constant value) was recently introduced and discussed in \cite{BV_2010, BV_2012}, where a classical inhomogeneous version of Kleinbock \& Margulis result was proven. Moreover,  Badziahin, Beresnevich  \& Velani \cite{BBV}  were later able to establish a general result implying the  following functional inhomogeneous statement for arbitrary non-degenerate manifolds. It contains both Theorems \ref{thm3} \& \ref{thm3a} as special cases ($\lambda \equiv  0$) - see also Theorem~1 of \cite{DB_inh_dim}.

\begin{thm}[Badziahin-Beresnevich-Velani]\label{thm4}
Let $\MMM \subset \RR^n$ be a non-degenerate manifold of dimension $m$ and let $\psi$ be an approximating function. Then, for any inhomogeneous function $\lambda$ satisfying $\lambda|_{\MMM} \in C^{(2)}$ and any $s \in (m-1,m]$ we have
$$
\mcal H^s(\W\cap \MMM) \: =\: \mcal H^s(\MMM) \quad {\rm if } \quad \sum\limits_{q=1}^{\infty}\psi(q)^{s-m+1}q^{n+m-s-1}=\infty.
$$
%
%
Furthermore, in the case that $s=m$ and so $\mathcal H^s \equiv \Lambda$ we have
$$
\mcal H^m(\W\cap\MMM)\:=\: 0 \quad\qquad {\rm if } \quad \sum\limits_{q=1}^{\infty}\psi(q)q^{n-1}<\infty.
$$
\end{thm}

At the time of writing, this represented the state of the art in solving the functional inhomogeneous analogue of the Generalized Baker-Schmidt Problem for Hausdorff measure. To reiterate, the following problem represents the main remaining hurdle in establishing a complete Hausdorff measure theory for Diophantine approximation on differentiable manifolds.

\begin{pro} \label{prob1}
Let $\MMM \subset \RR^n$ be a non-degenerate manifold of dimension $m$ and let $\psi$ be an approximating function. Then, for any inhomogeneous function $\lambda$ satisfying $\lambda|_{\MMM} \in C^{(2)}$ and any $s \in (m-1,m]$ verify that
$$
\mcal H^s(\W\cap \MMM) \: = \:  0 \quad {\rm if } \quad \sum\limits_{q=1}^{\infty}\psi(q)^{s-m+1}q^{n+m-s-1}<\infty.
$$
\end{pro}

\subsubsection{Further remarks}

Similar breakthroughs to those outlined above have been made in the setting of simultaneous Diophantine approximation on manifolds. Here, the focus is concentrated on the set $\mcal S(\lambda, \psi)$ introduced in \eqref{vb2} and progress has largely followed a parallel, although often slower moving, arc. For brevity we do not include full details of every result here, suffice to say that a complete analogue of Theorem~\ref{thm3a} and an analogue of the divergence part of Theorem~\ref{thm3} for simultaneous Diophantine approximation on analytic manifolds were both proven by Beresnevic in the recent seminal paper \cite{Ber_Ann}. This followed the establishing of these statements in the special case of  non-degenerate planar curves in both the case of convergence~\cite{VV} (for $C^{(2)}$ curves) and divergence~\cite{BDV_main} (for $C^{(3)}$ curves) - see also \cite{BZ}.  Classical inhomogeneous versions of the latter two results appeared shortly afterwards in the paper~\cite{BVV}. Consequently, the theory of simultaneous Diophantine approximation for $\mcal S(\lambda_0, \psi)\cap\mcal C$ is essentially complete. A simultaneous analogue of the convergence part of Theorem~\ref{thm3} and a simultaneous analogue of Theorem~\ref{thm4} for arbitrary non-degenerate analytic manifolds remain open problems.



We end our discussion with brief mention of a related research direction. In~1970, Baker $\&$ Schmidt {\cite{bs}} obtained a lower bound for the Hausdorff dimension of sets arising from Mahler's problem. They also conjectured that their bound was sharp. The original Baker-Schmidt conjecture was settled by Bernik~\cite{Bern83}. It was the extrapolation of this problem which led to the question of determining the Hausdorff measure of Diophantine sets restricted to a manifold. However, the special case of determining  the Hausdorff dimension of such sets has also been of great interest. We refer the reader to \cite{DB_inh_dim} and \cite{BBV} and the references therein for details and further discussion.


\subsection{Statements} Before stating our main result, we elaborate on the natural restrictions which we will impose on each planar curve.

Without loss of generality, we may assume that each $\CCC^{(2)}$ planar curve $\mcal C := \mcal{C}_f=\{(x, f(x)): \; x\in I\}$ can be given as
the graph of a  $C^{(2)}$ map $f : I\to \mb R,$ where $I$ is some compact interval of $\mb R$.  Given $s\in (0, 1]$ we say that a $C^{(2)}$ planar curve $\mathcal C$ is \textit{non-degenerate (with respect to $\mcal H^s$)} if the set of points on $\mathcal C$ at which the curvature vanishes is a set of $s$-dimensional Hausdorff measure zero; i.e., if $\mcal H^s\left(\{x\in I: f^{\prime\prime}(x)=0\}\right)=0$. It is readily verified that  in the case `$s=1$' this condition is equivalent to that required by Schmidt as described earlier - see \cite {VB, KM} for further discussion.

For $\CCC \in C^{(2)}(I)$, the above observations allow us to re-express the intersection $\W \cap \CCC$ as the one parameter set $W(\mcal C, \la, \psi)$ given by
\begin{equation*}
  W(\mcal C, \la, \psi)=\left\{x\in I: \begin{array}{l}
  |a_2f(x)+a_1x+a_0+\la|_\CCC(x)|<\psi(|\mbf a|) \\
  \text{for} \  i.m. \  (a_0, \mbf a) : =(a_0, a_1, a_2)\in \mb Z\times \mb Z^2\setminus\{\mbf 0\}
                                  \end{array}\right\}.
\end{equation*}
We may now state our main result.

\begin{thm}\label{thm6}
  Let $\psi$ be an approximating function and let $\CCC$ be a $C^{(2)}(I)$ non-degenerate planar curve. Then for any inhomogeneous function $\lambda$ satisfying $\lambda|_{\CCC} \in C^{(2)}$ and any $s \in (0, 1)$ we have
\begin{equation}\label{conresult}\mcal H^s( \w)=
  0 \ \  {\rm if } \quad \sum\limits_{q=1}^{\infty}\psi(q)^sq^{2-s}<\infty. \end{equation}
\end{thm}
 Together with Theorem~\ref{thm4}, we obtain the following combined corollary upon recalling that Hausdorff measure is invariant (up to a constant) under bi-Lipshitz bijections such as the co-ordinate map $\mbf f: I \to \CCC_f: x \rightarrow (x, f(x))$.
 \begin{cor} \label{cor}
    Let $\psi$ be an approximating function and let $\CCC$ be a $C^{(2)}(I)$ non-degenerate planar curve. Then for any inhomogeneous function $\lambda$ satisfying $\lambda|_{\CCC} \in C^{(2)}$ and any $s \in (0, 1]$ we have
\begin{equation*}
  \mcal H^s( \W \cap \CCC)=\left\{\begin{array}{cl}
 0 &  {\rm if } \quad\sum\limits_{q=1}^{\infty}\psi(q)^sq^{2-s} < \infty.\\[3ex]
  \mcal H^s(\mcal C) &  {\rm if } \quad \sum\limits_{q=1}^{\infty}\psi(q)^sq^{2-s}=\infty.
                                     \end{array}\right.
\end{equation*}
\end{cor}

\begin{rem}We remark that the assumption $s\in (0, 1)$ in Theorem \ref{thm6} is absolutely necessary. The case `$s=1$' in Corollary \ref{cor} corresponds to the result of Badziahin, Beresnevich \& Velani, and in the case `$s>1$' the string of inequalities  $\mcal H^s( \w)\leq \mcal H^s( \mcal C)=0$ holds for any approximating function~$\psi$ without the requirement of any sum condition.
\end{rem}

The above theorem brings the theory of dual approximation on non-degenerate curves in line with the advances in simultaneous Diophantine approximation on non-degenerate planer curves in~\cite{BVV}, and also generalises the result of Huang~\cite{Huang} concerning dual Diophantine approximation on non-degenerate planer curves to the functional inhomogeneous setting. In doing so, we provide a complete solution to Problem~1 in the special case `$\MMM=\CCC$, $n=2$'. Indeed, our proof synthesises the ideas of~\cite{Huang} with those presented in~\cite{DB_inh_dim}.


\section{Preliminaries}\label{Hausdorffdefs}

To simplify notation in the proofs the Vinogradov symbols $\ll$ and $\gg$ will be used to indicate
 an inequality with an unspecified positive multiplicative constant.
 If $a\ll b$ and $a\gg b$  we write $a\asymp
 b$, and say that the quantities $a$ and $b$ are comparable.

\subsection{Hausdorff measure and dimension}\label{hm}

For completeness we give a very brief introduction to Hausdorff measures and dimension. For further details see \cite{berdod}.
Let
$F\subset \mb{R}^n$.
 Then, for any $\rho>0$ a countable collection $\{B_i\}$ of balls in
$\mb{R}^n$ with diameters $\mathrm{diam} (B_i)\le \rho$ such that
$F\subset \bigcup_i B_i$ is called a $\rho$-cover of $F$.
Let
\[
\mcal{H}_\rho^s(F)=\inf \sum_i \mathrm{diam}(B_i)^s,
\]
where the infimum is taken over all possible $\rho$-covers $\{B_i\}$ of $F$. It is easy to see that $\mcal{H}_\rho^s(F)$ increases as $\rho$ decreases and so approaches a limit as $\rho \rightarrow 0$. This limit could be zero or infinity, or take a finite positive value. Accordingly, the \textit{Hausdorff $s$-measure $\mcal{H}^s$} of $F$ is defined to be
\[
\mcal{H}^s(F)=\lim_{\rho\to 0}\mcal{H}_\rho^s(F).
\]
It is easily verified that Hausdorff measure is monotonic and countably sub-additive, and that $\mcal{H}^s(\emptyset)=0$. Thus it is an outer measure on $\RR^n$. Furthermore, for any subset $F$ one can easily verify that there exists a unique critical value of $s$ at which $\mcal{H}^s(F)$ `jumps' from infinity to zero. The value taken by $s$ at this discontinuity is referred to as the \textit{Hausdorff dimension of  $F$} is denoted by $\dim F $; i.e.,
\[
\dim F :=\inf\{s\in \mb{R}^+\;:\; \mcal{H}^s(F)=0\}.
\]

When s is an integer, $n$ say, then $\mcal{H}^n$ coincides with standard $n$-dimensional Lebesgue measure, which we will denote $|\, . \, |$. 
Hausorff $s$-measure, like Lebesgue measure, is preserved (up to a constant) by certain well behaved maps.  In particular, if $g: E \to F$ is a bi-Lipshitz bijection between two sets in Euclidean space then $\HHH^s(E) \asymp \HHH^s(F)$.

\subsection{Auxiliary lemmas} We now group together two important results that we appeal to in the course of proving Theorem \ref{thm6}.
Firstly,  we state the Hausdorff measure version of the famous Borel-Cantelli lemma (see Lemma 3.10 of \cite{berdod}) which will be key to our method. It allows us to estimate the Hausdorff measure of certain sets via calculating the $s$-Hausdorff sum of a `nice' cover.
\begin{lem}\label{bclem}
    Let $H_i$ be a sequence of hypercubes in $\mb R^n$ and suppose that for some $s>0$ we have $$\sum_i |H_i|^s \, < \, \infty.$$ Then, $$\mcal H^s\left(\bigcap_{r=1}^\infty \bigcup_{i=r}^\infty H_i\right)=0.$$
\end{lem}

Secondly, the following famous lemma can be found in \cite[Lemma
 9.7]{Harman} and is originally attributed to Pyartly. We will use it several times throughout the proof.
\begin{lem}[\cite{Harman}]\label{Pya_lem} Let $h(x)\in C^{(2)}(I)$ be such that $\min\limits_{x\in I}|h^\prime(x)|=\delta_1$ and $\min\limits_{x\in I}|h^{\prime\prime}(x)|=\delta_2$. For $\eta>0$, define
\begin{equation*}
  E(\eta):=\{\eta\in I: |h(x)|<\eta\}.
\end{equation*}
Then we have
\begin{equation*}
  |E(\eta)|\ll \min\left(\frac\eta\delta_1, \sqrt{\frac\eta\delta_2}\right).
\end{equation*}
\end{lem}


\subsection{Further observations}

We next make some standard simplifications. Without loss of
generality we can assume that there exist positive absolute
constants $c_1$ and $c_2$ such that
\begin{equation}\label{conv}
c_1\leq |f^{\prime\prime}(x)| \leq c_2 \qquad \forall x\in I.
\end{equation}
Indeed, since the set $\{x\in I\;:\; f''(x) = 0\}$ is compact, its
complement $\{x\in I: f^{\prime\prime}(x)\neq0\} $ is a countable
union of intervals $J_\alpha, \alpha\in \mathbb{Z}$ such that
$$
\forall \alpha\in \mathbb{Z}\;\; \exists
c_1(\alpha),c_2(\alpha)>0\;\mbox{ such that }\; \forall x\in
J_\alpha,\;c_1(\alpha)\le  |f''(x)|\le c_2(\alpha).
$$
So if we prove that $\mcal H^s(W(\mcal C,\lambda,\psi))=0$ on every
interval $J_\alpha$ then this together with the fact that $\mcal
H^s\left(\{x\in I: f^{\prime\prime}(x)=0\}\right)=0$ will give us
that $W(\mcal C,\lambda,\psi)$ has zero Hausdorff $s$-measure on the
whole interval $I$. Additionally, since $\lambda|_\CCC$ is $C^{(2)}$ we may also assume a similar set of inequalities hold for this restriction of the inhomogeneous function. For convenience, we will from here on abuse notation by simply writing $\lambda$ for  $\lambda|_\CCC$.

%

As required by the theorem, we assume throughout that
\begin{equation}\label{coneq}
\sum\limits_{q=1}^{\infty}\psi(q)^sq^{2-s}<\infty.
\end{equation}
Since $\psi$ is a monotonic function, a Cauchy condensation argument yields that
\begin{equation}\label{condensation}
\sum\limits_{t=1}^{\infty}\psi(2^t)^s \: 2^{(3-s)t} \: \: \asymp \: \: \sum\limits_{q=1}^{\infty}\psi(q)^sq^{2-s} \: \: < \: \: \infty.
\end{equation}
In particular, we have
$\psi(2^t)^s\: 2^{(3-s)t}\to0$, whence
\begin{equation}\label{vb6}
\psi(q)< q^{1-\frac 3 s}\qquad\text{for sufficiently large }q.
\end{equation}
On the other hand, by following the arguments from \cite{Huang},
without loss of generality we may assume for small
$\epsilon>0$ that
\begin{equation}\label{psi_low}
\psi(q)\ge q^{1-\frac{3+\epsilon}{s}}.
\end{equation}
Indeed, if $\psi$ does not satisfy \eqref{psi_low} one can instead  consider the
auxiliary function $$\Psi(q):=\max\left\{\psi(q),
q^{1-\frac{3+\epsilon}{s}} \right\},$$ which is an approximating
function satisfying both \eqref{psi_low} and $$
\sum\limits_{q=1}^{\infty}\Psi(q)^sq^{2-s}<\infty.$$

\section{Proof of Theorem~\ref{thm6}}

Our general strategy will be to carefully construct a cover of each set $W(\mcal C, \la, \psi)$ with a collection of bounded intervals. We will then estimate the measure of each element of the cover and their number.  This measure will not be uniform over all types of interval in the cover and the number of intervals of a given size will vary. Calculating an estimate for every type of interval will then allow us to find an estimate for the Hausdorff measure of the entire set. 

For notational convenience let $F(x):=F_{(a_0, \mbf a)}(x):=a_2f(x)+a_1x+a_0$ for any fixed triple $(a_0, \mbf a) : =(a_0, a_1, a_2)\in \mb Z\times \mb Z^2\setminus\{\mbf 0\}$ . Given any such function $F$ denote
\begin{equation*}
\Delta(F, \la, \psi):=\{x\in I: |F(x)+\la(x)|<\psi(|\mbf a|)\}.
\end{equation*}
It follows that the sets $\Delta(F, \la, \psi)$ form a cover of $\w$; indeed, we have that
\begin{equation}\label{limsup1}
  \w \: = \: \bigcap_{N \in \NN} \: \bigcup_{r \geq N} \: \bigcup_{\substack{ (a_0, \mbf a): \\  |\mbf a|=r}}{}{}\Delta(F, \la, \psi),
\end{equation}
is precisely the set of points $x \in I$ that lie in infinitely many of the sets $\Delta(F, \la, \psi)$. However, to apply Lemma~\ref{bclem} we must find a cover by intervals (one-dimensional hypercubes). One can readily verify that each set $\Delta(F, \la, \psi)$ comprises a finite union of disjoint intervals. Given a fixed triple $(a_0, \mbf a)$ we denote by $K(F):=K(F, \lambda, \psi)$ the collection of such intervals~$H$ inside $\Delta(F, \la, \psi)$ and  by $k(F)$ their number. In other words, we have $$\Delta(F, \la, \psi) = \bigsqcup_{H \, \in \, K(F)} \, H \quad \text{ and } \quad k(F)=\#K(F).$$ In this way, a cover for the set $\w$ by intervals has easily been constructed for each $N \in \NN$. Let $\mathcal K(N)$ be the collection of all intervals in this cover; that is, let
$$
\mathcal K(N) \, = \, \bigcup_{r \geq N}\bigcup_{\substack{ (a_0, \mbf a):  \\ |\mbf a|=r}} \: \:K(F)$$
By Lemma \ref{bclem} and reformulation~\eqref{limsup1}, in order to prove Theorem~\ref{thm6} it therefore suffices to demonstrate that for some $N\in \NN$ we have
\begin{equation} \label{endgame}
\ws: \: = \: \sum_{H \, \in \, \mathcal K(N)} \, |H|^s \: \ll \: \sum_{r \geq N}\sum_{\substack{ (a_0, \mbf a):  \\ |\mbf a|=r}} \: \:\sum_{H \, \in \, K(F)} |H|^s \: \:  < \: \: \infty.
\end{equation}

Notice that for any fixed  $(a_0, \mbf a)$ we have
$$
\sum_{H \, \in \, K(F)} |H|^s    \: \leq k(F) \cdot |\Delta(F, \la, \psi)|^s.
$$
Moreover, as we will see later the number of intervals $k(F)$ for any $F$ is uniformly bounded by an absolute constant. This implies that for any collection $I_0$ of triples~$(a_0, \mbf a)$ we have
\begin{equation}\label{hypercubecover}
\sum_{ (a_0, \mbf a) \, \in \, I_0} \: \: \sum_{H \, \in \, K(F)} |H|^s \: \ll \: \sum_{ (a_0, \mbf a) \, \in \, I_0}|\Delta(F, \la, \psi)|^s.
\end{equation}
In many cases this crude estimate will be enough for us, but there are many instances where we will first need to refine the cover formed by the intervals in $\mathcal{K}(N)$ by refining certain collections  $K(F)$.
During such a process it should be understood that demonstrating the final inequality of \eqref{endgame} holds for any refinement of $\mathcal{K}(N)$ will still be sufficient to prove Theorem~\ref{thm6}.

In practice  we will only be interested in polynomials $F$ with $\Delta( F,
\la, \psi)\neq\emptyset$, so certainly we may assume for $\mbf a= (a_1, a_2)$ large enough that  $|F_{(a_0, \mbf a)}(x)+\la (x)|<\psi(|\mbf
a|)< |\mbf a|^{1-3/s} <1$ for some $x\in I$. Consequently,  for a fixed $\mbf a$
there are only finitely many~$a_0$ such that $\Delta(F_{(a_0, \mbf a)}(x), \la, \psi)\neq
\emptyset.$ Moreover, upon setting $$\kappa:=\underset{x\in \emph
I}{\max}\left\{|x|+|f(x)|+|\la(x)|+1\right\}$$ we have for every such $a_0$ that
\begin{equation}\label{a0counting}
|a_0|\leq \kappa |\mbf a|.\end{equation}

Now, for every fixed $\mbf a$ we may assume
\begin{equation*}
  |F^\prime(x)+\la^\prime(x)|\geq |a_1|-|a_2f^\prime(x)+\la^\prime(x)|\geq |a_1|-2|a_2|M,
\end{equation*}
where $M:=\underset{x\in \emph  I}{\max}\left\{|f^\prime(x)|+|\la^\prime(x)|\right\}$ is a fixed constant.
As a consequence, we split the possible choices for the pair of integers $(a_1, a_2)\in\mb Z^2\setminus \left\{ \mbf 0 \right\}$ into two exhaustive subsets; let
\begin{equation*}
  I_1:=\{(a_1, a_2)\in\mb Z^2\setminus \left\{ \mbf 0 \right\}: |a_1|>4M|a_2|\},
\end{equation*} and
\begin{equation*}
  I_2:=\{(a_1, a_2)\in\mb Z^2\setminus \left\{ \mbf 0 \right\}: |a_1|\leq 4M|a_2|\}.
\end{equation*}
One can see that for each $(a_1,a_2) \in I_1$ and for each $x\in I$ we have
\begin{equation}\label{I1}
|F'(x)+\lambda'(x)|\ge |a_1|/2.
\end{equation}
Futhermore,
\begin{equation}\label{I1a}
 |\mbf a|\asymp|a_1|
\end{equation}
for each $\mbf a=(a_1,a_2)$ in $I_1$.

We now analyze the sum $\ws$ arising in \eqref{endgame} and split the summand into two natural cases according to whether $\mbf a$ lies in $I_1$ or $I_2$. To be precise, let
\begin{equation*}\label{rr}
\wrp=\sum_{r \geq N}{} \: \sum_{\substack{ (a_0, \mbf a): \\ \mbf a \, \in \, I_1, \: |\mbf a|=r}}|\Delta(F, \la, \psi)|^s,
\end{equation*}
and $$\wds=\sum_{r \geq N}{} \: \sum_{\substack{ (a_0, \mbf a): \\ \mbf a \, \in \, I_2, \: |\mbf a|=r}}\sum_{H \, \in \, K(F)} |H|^s .$$
In view of \eqref{hypercubecover} one can readily verify that
\begin{equation}\label{mh1}
\ws  \: \ll \: \wrp \: + \: \wds,
\end{equation}
whence, the desired statement that $\mcal H^s(\w)=0$ will follow upon
establishing separately
\begin{center}
  \[ \mbf{Case \ I:}\qquad \qquad \wrp < \infty.\]
  \[\mbf{Case \ II:} \ \ \quad \qquad \wds < \infty.\]
\end{center}
We now split our proof into two parts, each corresponding to one of the above cases.

\subsection{Establishing Case I}

%
We begin by establishing the  much easier Case I.  In doing so we appeal to Lemma~\ref{Pya_lem}, one consequence of which is that for every $\mbf a=(a_1,a_2) \in I_1$ we have
\begin{equation*}
|\Delta(F, \la, \psi)|\ll  \frac{\psi(|\mbf a|)}
{\min_{x \in I}|f'(x)+\lambda'(x) |   } \: \stackrel{\eqref{I1}}{\leq}  \:     \frac{2\psi(|\mbf a|)}
{|a_1|}\:\stackrel{\eqref{I1a}}{\asymp} \:\frac{\psi(|\mbf a|)}{|\mbf a|}.
\end{equation*}
For sufficiently large $N$ it follows from \eqref{a0counting} and \eqref{I1a} that there are at most $\ll r^2$ triples $(a_0, a_1, a_2)$ with $\mbf a \in I_1$ and $|\mbf{a}| = r$.
Furthermore, it follows from \eqref{I1} and $\eqref{conv}$ that each set $\Delta(F, \la, \psi)$ with $\mbf a \in I_1$ comprises a single interval; i.e., $k(F) = 1$. Hence,
\begin{eqnarray*}\label{fvol2}
\wrp \: \ll \:  \sum_{r\geq N}\sum_{\substack{(a_0, \mbf a): \\ \mbf a \in I_1, \:|\mbf a|=r}}\left(\frac{\psi(|\mbf a|)}{\mbf a}\right)^s & = &
\sum_{r\geq N}\left(\frac{\psi(r)}{r}\right)^s \sum_{\substack{(a_0, \mbf a): \\ \mbf a \in I_1, \:|\mbf a|=r}} \, 1
\\[1ex]
 &\ll& \sum_{r\geq N}\psi(r)^s r^{-s} \cdot r^2<\infty,
\end{eqnarray*}
by \eqref{coneq}
as required.

\subsection{Preliminaries for Case II}
 Establishing Case II is much more difficult and will involve various further refinements to the cover. Firstly, observe that in Case II we have by the definition of the set $I_2$ that
\begin{equation}\label{I2a}
a_2\asymp |\mbf a|.
\end{equation}
For notational convenience we denote   $\F(x):=\F_{\mbf a}(x):= a_1 x + a_2f(x) + \lambda(x)$,  so that for any $a_0 \in \ZZ$ we have $\F_{\mbf a}(x) + a_0 = F_{(a_0, \mbf a)}(x)$. Since $\lambda''$ is bounded on $I$ from above and below by absolute constants (that is, $|\lambda''(x)|\ll 1$ for $x \in I$),  inequalities~\eqref{conv} and~\eqref{I2a} imply that for $|\mbf a|$ large enough
\begin{equation}\label{observation}
 |\F''(x)|=|a_2f''(x) + \lambda''(x)|\asymp |\mbf a|  \qquad \text{ for } x\in I.
\end{equation}
This observation implies that each set $\Delta(F, \la, \psi)$ arising from a triple $(a_0, \mbf a)$ with $\mbf a \in I_2$ large enough comprises at most two intervals; i.e.,  $k(F) \leq 2$ for all such $(a_0, \mbf a)$. It also follows that for $|\mbf a|$ large enough the function $\F'(x)=a_2f'(x)+a_1 + \lambda'(x)$ is monotonic. 
Given a pair $(a_1,a_2) \in I_2$ we
define $x_0:=x_0(\mbf a)$ as the root of the equation
$$
\F^\prime(x) = 0.
$$
If such a root does not exist we conventionally set $x_0=\infty$ and $\F(x_0)=\infty$.

As in \cite{Huang}, we proceed by assuming $a^\prime_0$ to be a unique integer such that that $-1/2 < \F(x_0)
- a^\prime_0\le 1/2$. For any given triple $(a_0, a_1, a_2)$ there are two possibilities; either $a_0\neq a^\prime_0$ or $a_0= a^\prime_0$. The final case will prove to be the most difficult scenario to deal with therefore we split our calculation into two further subcases. In fact, we will repeat this procedure several more times in the latter subcase. Observe that \eqref{hypercubecover} implies that
\begin{equation}
\wds \: \ll \: \wdsa \: + \: \wdsb,
\end{equation}
where
\begin{equation*}\label{rr2}
\wdsa=\sum_{r \geq N}{} \: \sum_{\substack{ (a_0, \mbf a): \: a_0\neq a^\prime_0,\\ \mbf a \, \in \, I_2, \: |\mbf a|=r}}|\Delta(F, \la, \psi)|^s,
\end{equation*}
and $$\wdsb=\sum_{r \geq N}{} \: \sum_{\substack{ (a_0, \mbf a): \: a_0= a^\prime_0,\\ \mbf a \, \in \, I_2, \: |\mbf a|=r}} \: \: \sum_{H \, \in \, K(F)} |H|^s.
$$
Thus,
the desired statement that $\wds<\infty$ will follow upon
establishing separately
\begin{center}
  \[ \mbf{Case \ IIa:}\qquad \qquad \wdsa < \infty.\]
  \[\mbf{Case \ IIb:} \ \ \quad \qquad \wdsb < \infty.\]
\end{center}
As before, we split the remaining proof into two distinct parts.

\subsection{Establishing Case IIa} We begin with an important lemma.
\begin{lem}\label{lem2}
For each $x\in I$ and for every $\mbf a=(a_1, a_2) \in I_2$ with $|\mbf a|$ sufficiently large we have
$$
|\F'(x)| \asymp |a_2(\F(x) - \F(x_0))|^{1/2}
$$
and
$$|\F(x)-\F(x_0)| \: \gg \: |a_0-a'_0|.$$
\end{lem}
A proof of this lemma follows from the proof of~\cite[Lemma  3]{Huang} and the calculation immediately succeeding it with straightforward modifications. For brevity we do not include full details here,
%
suffice to say that one only needs to incorporate observation~\eqref{observation}.



%

 By applying Lemma \ref{lem2} we have for any $(a_0, \mbf a) \in \ZZ \times I_2$ with $|\mbf a|$ sufficiently large that
$$ |\F'(x)| \: \gg \: (|a_2| \cdot|a_0-a'_0|)^{1/2}.$$
In turn, Lemma \ref{Pya_lem} then yields
\begin{equation}\label{I2}
 |\Delta(F, \la, \psi)| \: \ll   \:     \frac{\psi(|\mbf a|)}
 {(|a_2| \cdot |a_0-a'_0|)^{1/2}}  \: \stackrel{\eqref{I2a}}{\asymp}   \: \frac{\psi(|\mbf a|)}
 {(|\mbf a| \cdot |a_0-a'_0|)^{1/2}}.
\end{equation}

Next, for a given $r \geq N$  there are at most $\ll r$  
choices for $(a_1, a_2)\in I_2$ with $|\mbf a| = r$. Hence, for $N$ sufficiently large it follows from \eqref{a0counting} and \eqref{I2}  that

\begin{eqnarray*}
  \wdsa&\ll& \sum_{r\geq N}\: \:\sum_{\substack{(a_0, \mbf a): \: a_0\neq a^\prime_0,\\ \mbf a \: \in \: I_2, \:|\mbf a|=r }}
|\Delta (F, \la, \psi)|^s\\ &\ll&  \sum_{r\geq N}\: \:\sum_{\substack{\mbf a\in I_2: \\ |\mbf a|=r}} \: \sum_{\substack{a_0 \leq \kappa r: \\a_0\neq a^\prime_0}}
\left(\frac{\psi (|\mbf a|)}{|\mbf a|^{1/2}}\left(|a_0-a_0^\prime|\right)^{-1/2}\right)^s\\ &=&\sum_{r\geq N} \left(\frac{\psi (r)}{r^{\frac12}}\right)^s \: \: \sum_{\substack{a_0 \leq \kappa r: \\a_0\neq a^\prime_0}}\left(|a_0-a_0^\prime|\right)^{-s/2}  \sum_{\substack{(a_1,a_2)\in I_2: \\ |\mbf a|=r}} \:  \: 1\\ &\ll& \sum_{r\geq N} \psi (r)^s \, r^{1-\frac{s}{2}} \sum_{\substack{a_0 \leq \kappa r: \\a_0\neq a^\prime_0}}\left(|a_0-a_0^\prime|\right)^{-s/2}. \
\end{eqnarray*}
Finally, observe that for each $r \geq N$ the final sum
$$
\sum_{\substack{a_0 \leq \kappa r \\a_0\neq a^\prime_0}}\left(|a_0-a_0^\prime|\right)^{-s/2} \: \ll \: \sum_{t=1}^{\kappa r} t^{-s/2},
$$
which is  by partial summation bounded above by a constant times $r^{1-s/2}$ for sufficiently large $N$. Hence, as required the above calculation yields
$$
 \wdsa \:  \ll \:  \sum_{r\geq N} \psi (r)^s \, r^{1-\frac{s}{2}} \cdot r^{1-s/2} = \sum_{r\geq N} r^{2-s}\psi (r)^s \, \stackrel{\eqref{coneq}}{<} \, \infty.
$$

\subsection{Preliminaries for Case IIb} The proof in Theorem \ref{thm6} now follows upon finding a suitable bound for the sum
\begin{equation}\label{finalone}
\wdsb=\sum_{r \geq N}{} \: \sum_{\substack{ (a_0, \mbf a): \: a_0= a^\prime_0,\\ \mbf a \, \in \, I_2, \: |\mbf a|=r}} \: \:\sum_{H \, \in \, K(F)} |H|^s.
\end{equation}
For notational convenience, for $r \geq N$ let $A_r$ denote the set of triples  appearing in the above sum with $|\mbf a|=r$; that is, let
$$
A_r:\: = \: \left\{(a_0, \mbf a) \in \ZZ \times I_2: \: a_0= a^\prime_0, |\mbf a|=r \right\}.
$$
Note that there is only one choice for $a_0$ for a given $\mbf a \in I_2$. In other words, for any $r \geq N$ we have $\#A_r \ll r$.

Case IIb is by far the most difficult and will require rather intricate treatment. 
We will first need to decompose the collections $K(F)$ in a  sophisticated way.  
How we deal with points in a set $\Delta(F,\lambda, \psi)$ for some triple $(a_0, \mbf a)$ appearing in the summand of~$(\ref{finalone})$ will depend upon the behaviour of the derivative $\F'(x)$ on each of the intervals $H\in \, K(F)$ of which  $\Delta(F,\lambda, \psi)$ comprises. For this reason we proceed by first dividing the collection of sets $\Delta(F,\lambda, \psi)$ associated with $(\ref{finalone})$ into two exhaustive subsets. For each triple $(a_0, \mbf a)$ under consideration let
$$
\Delta_1(F,\lambda, \psi):= \{x\in \Delta(F,\lambda, \psi)\;:\; |\F'(x)|\ge |\mbf
a|^{1-1/s}\}
$$
and $$\Delta_2(F,\lambda, \psi): = \Delta(F,\lambda, \psi)\backslash
\Delta_1(F,\lambda).$$
This process will inevitably split points in some intervals $H\in \, K(F)$ between the two categories. To be specific, in view of observation \eqref{observation} the set  $\Delta_1(F,\lambda, \psi)$ naturally comprises of a collection $K_1(F)$ of at most four intervals of the form $H \cap \Delta_1(F,\lambda, \psi)$  for some $H\in \, K(F)$. For this reason, as in \eqref{hypercubecover} we will in this case freely use the estimate
\begin{equation}\label{hypercubecover1}
\sum_{ (a_0, \mbf a) \, \in \, A_r} \: \: \sum_{H \, \in \, K_1(F)} |H|^s \: \ll \: \sum_{ (a_0, \mbf a) \, \in \, A_r}|\Delta_1(F, \la, \psi)|^s.
\end{equation}
It is also the case that each set $\Delta_2(F,\lambda, \psi)$ may consist of at most four intervals of the form $H \cap \Delta_2(F,\lambda, \psi)$, the collection of which we denote by $K_2(F)$. However, an estimate such as \eqref{hypercubecover1} will prove too crude to apply in our method. We will first need to further refine the cover.

Without loss of generality we may refine the cover formed by intervals in the collection $\mathcal K(N)$ by replacing intervals from a given $K(F)$ by intervals from $K_1(F)$ and $K_2(F)$ in the volume sum $\wdsbb$.  Accordingly,  we may specialise once again by setting
\begin{equation*}\wdsba=\sum_{r\ge N}\: \sum_{(a_0, \mbf a) \: \in \: A_r}
|\Delta_1(F,\lambda, \psi)|^s,
\end{equation*}
and
\begin{equation*}\wdsbb=\sum_{r\ge N}\: \sum_{(a_0, \mbf a) \: \in \: A_r}
\: \:\sum_{H \, \in \, K_2(F)} |H|^s.
\end{equation*}
It follows from observation \eqref{hypercubecover1} that
$$\wdsb \: \ll \: \wdsba \: + \: \wdsbb.$$ We now aim to show that the two latter series converge.

\subsubsection{Establishing subcase IIb1}
Since we are assuming $|\F'(x)|\ge |\mbf a|^{1-1/s}$ for every $x \in \Delta_1(F,\lambda, \psi)$ it follows from Lemma~\ref{Pya_lem} that
$$
|\Delta_1(F,\lambda, \psi)| \: \ll \: \frac{\psi(|\mbf a|)}{|\mbf a|^{1-1/s}}.
$$
Recalling that for any $r \geq N$ the number of triples $(a_0, a_1, a_2) \in \ZZ \times I_2$ for which $|\mbf a|=r$ is at most $\ll r$, this in turn implies that for $N$ sufficiently large we have
\begin{eqnarray*}
  \wdsba&\ll& \sum_{r\geq N}\: \: \sum_{(a_0, \mbf a) \: \in \: A_r}
|\Delta_1(F, \la, \psi)|^s\\ & \ll &  \sum_{r\geq N} \: \:\sum_{(a_0, \mbf a) \: \in \: A_r}
\left(\frac{\psi (|\mbf a|)}{|\mbf a|^{1-1/s}}\right)^s \\ &\ll& \sum_{r\geq N} \psi(r)^s \, r^{2-s}\: < \: \infty.
\end{eqnarray*}

\subsection{Establishing subcase IIb2}

All that remains is to demonstrate the convergence of the sum~$\wdsbb$. In turn this will establish Case IIb, Case II and Theorem~\ref{thm6} respectively. In this subcase  we have $|\F'(x)|< |\mbf a|^{1-1/s}$. Recall that observation~\eqref{observation} implies that we may assume $|\F''(x)|\gg |\mbf a|$ for $|\mbf a|$ sufficiently large. So, by  Lemma~\ref{Pya_lem} we have
\begin{equation}\label{dflam}
|\Delta_2(F,\lambda, \psi)| \ll \left(\frac{\psi(|\mbf a|)}{|\mbf
a|}\right)^{1/2}.
\end{equation}
First, we split  the triples into dyadic blocks $2^t\le |\mbf
a|<2^{t+1}$ and consider all functions~$\F$  arising in this way. For notational convenience, let
$
A^t:= \bigcap_{r=2^t}^{2^{t+1}-1} \: A_r; 
$
that is, let
$$
A^t:\: = \: \left\{(a_0, \mbf a) \in \ZZ \times I_2: \: a_0= a^\prime_0, |\mbf a|=r, \, 2^t \leq r < 2^{t+1} \right\}.
$$
Then, we may re-express the sum $\wdsbb$ as
\begin{equation}\label{case_iv}
\wdsbb \: \asymp \: \sum_{t\, \geq \, \left\lfloor\log_2N \right\rfloor} \:\:\sum_{(a_0, \mbf a) \: \in \: A^t}
\: \:\sum_{H \, \in \, K_2(F)} |H|^s.
\end{equation}

Before proceeding, we must further refine the intervals from the collections $K_2(F)$ appearing in~\eqref{case_iv}. Let $c = 1+\epsilon_1$, where $\epsilon_1$ is some very small
positive constant. For every $t$ we split the interval $I$ into $\ll 2^{ct}$
subintervals $J_i^{(t)}$ of the same length. So, for each such subinterval $J_i^{(t)} \subset I$ we have
\begin{equation}\label{Jsize}
|J_i^{(t)}| \ll 2^{-ct}.
\end{equation}
We will denote by $J(t)$ the number of subintervals created (and so $J(t)\ll 2^{ct}$).

We refine $\mathcal{K}(N)$ by replacing  each collection of intervals $$K_2(F) \: \:=\: \: \left\{ H \cap \Delta_2(F,\lambda, \psi): \: H \, \in \, K(F)\right\}\: \: \subset\: \: \mathcal{K}(N)$$
 appearing in \eqref{case_iv} with the collections
$$
K^i_2(F, \, \la): \, = \,  \left\{ H^\prime \cap J_i^{(t)}: \: H^\prime \, \in \, K_2(F)\right\},
$$
for $i=1, \ldots, J(t)$. Redefining $\wdsbb$ accordingly, we have
\begin{equation*}
\wdsbb \quad \asymp \quad  \displaystyle \sum_{t\, \geq \, \left\lfloor\log_2N \right\rfloor}\: \: \sum_{i=1}^{J(t)}\: \:\sum_{(a_0, \mbf a) \: \in \: A^t}
\: \:\sum_{H \, \in \, K^i_2(F, \, \la)} |H|^s,
\end{equation*}
and in view of the discussion immediately succeeding \eqref{hypercubecover} it suffices to show that this new series converges. It will be sufficient to utilise the crude estimate $\#K^i_2(F, \, \la) \leq \#K_2(F) \leq 4$, yielding
\begin{equation}\label{yetanother}
\wdsbb \quad \ll \quad \displaystyle\sum_{t\, \geq \, \left\lfloor\log_2N \right\rfloor}  \: \: \sum_{i=1}^{J(t)}\: \:\sum_{(a_0, \mbf a) \: \in \:A^t}\: \: |\Delta_2(F,\lambda, \psi) \cap J_i^{(t)}|^s.
\end{equation}

Next, we fix another very small parameter $\epsilon_2>0$ and divide the subintervals~$J^{(t)}$ (for $i=1, \ldots J(t)$) into two categories in the following way: 

\begin{itemize}

\item \emph{ \textbf {`Good' Intervals}}. An interval $J^{(t)} \subset I$ is called \textit{Good} if  it intersects at most $2^{(3/2 - c - \epsilon_2)t}$ sets of the form $\Delta_2(F_{(a_0, \mbf a)},\lambda, \psi)$ for $(a_0, \mbf a) \in A^t$.
\item  \emph{ \textbf { `Bad' Intervals.}} An interval $J^{(t)} \subset I$ is called \textit{Bad} if  it intersects more than $2^{(3/2 - c - \epsilon_2)t}$ sets of the form $\Delta_2(F_{(a_0, \mbf a)},\lambda, \psi)$ for $(a_0, \mbf a) \in A^t$.
\end{itemize}
In accordance with this classification we split our volume sum once more into two parts. Indeed,  setting
\begin{equation*}
\wdsbbg= \displaystyle\sum_{t\, \geq \, \left\lfloor\log_2N \right\rfloor}  \: \: \sum_{\text{Good intervals } J^{(t)} }\: \:\sum_{(a_0, \mbf a) \: \in \:A^t}
\: \: |\Delta_2(F,\lambda, \psi) \cap J^{(t)}|^s,
\end{equation*}
and
\begin{equation*}\label{case_iv3}
\wdsbbb=  \displaystyle\sum_{t\, \geq \, \left\lfloor\log_2N \right\rfloor}  \: \: \sum_{\text{Bad intervals } J^{(t)} }\: \:\sum_{(a_0, \mbf a) \: \in \:A^t}
\: \:\: \: |\Delta_2(F,\lambda, \psi) \cap J^{(t)}|^s,
\end{equation*}
it is easy to see in view of (\ref{yetanother}) that
\begin{equation*}\label{case_iv1}
\wdsbb \: \ll \:  \wdsbbg \: + \: \wdsbbb.
\end{equation*}
Once more we split our proof into two parts.

\subsubsection{The case of Good intervals}

We first estimate the measure of the sets $\Delta_2(F, \la, \psi)$ lying inside a given Good interval $J^{(t)}$. Recall that $\psi$ is decreasing and for every $t$ the number
of sets $\Delta_2( F, \la, \psi)$ satisfying $(a_0, \mbf a) \in A^t$ and $\Delta_2( F, \la, \psi)\cap J^{(t)}\neq \emptyset$ is by definition
$\ll 2^{(3/2 - c - \epsilon_2)t} $. Also note that by~\eqref{dflam} we have for every $(a_0, \mbf a) \in A^t$ and every Good or Bad interval $J^{(t)}$  that
\begin{equation}\label{dflam1}
|\Delta_2(F,\lambda, \psi)\cap J^{(t)}|^s \: \: \ll\: \: \left(\frac{\psi(|\mbf a|)}{|
    \mbf a|}\right)^{s/2}\: \: \ll \: \: \psi(2^t)^{s/2} \: 2^{-st/2}.
\end{equation}
Recalling that the total number  $J(t)$ of possible intervals, Good or Bad, satisfies $J(t)\ll 2^{ct}$, we certainly have for $N$ sufficiently large that
\begin{eqnarray*}
\wdsbbg &=&\sum_{t\, \geq \, \left\lfloor\log_2N \right\rfloor}  \: \: \sum_{\text{Good intervals } J^{(t)} }\: \:\sum_{(a_0, \mbf a) \: \in \:A^t}
\: \: |\Delta_2(F,\lambda, \psi) \cap J^{(t)}|^s \\
 &\ll& \sum_{t\, \geq \, \left\lfloor\log_2N \right\rfloor} J(t)\cdot 2^{(3/2 -
c - \epsilon_2)t} \cdot \psi(2^t)^{s/2}\: 2^{-st/2}\\ &\ll&
 \sum_{t\, \geq \, \left\lfloor\log_2N \right\rfloor} \:\psi(2^t)^{\frac s2} \cdot   2^{(\frac3 2-\frac s 2-\epsilon_2)t}\:
\\ & \stackrel{\eqref{vb6}}{\ll}& \sum_{t\, \geq \, \left\lfloor\log_2N \right\rfloor} 2^{-\epsilon_2t} \quad \ll  \quad  \sum_{r\, \geq \,N} r^{-(1+\epsilon_2)} \quad < \quad \infty. \
\end{eqnarray*}


\subsubsection{The case for Bad intervals}

Assume we are given a Bad interval $J^{(t)}$ and some $(a_0, \mbf a) \in A^t$ for which $\Delta_2(F,\lambda, \psi)$ intersects $J^{(t)}$. Then, there must exist a point $y_0\in \Delta_2(F,\lambda, \psi) \cap J^{(t)}$ satisfying
\begin{equation}\label{bounds}
|\F(y_0)| <\psi(2^t) \quad \text{ and } \quad |\F'(y_0)| < 2^{(1-1/s)t}.
\end{equation}
Since $s \in (0,1)$, for $\epsilon_1$ small
enough we may assume that $3/s-1 > 2+3\epsilon_1$. This and~\eqref{vb6} together imply that 
\begin{equation}\label{psibound}
\psi(2^t) < 2^{(1-3/s)t} <  2^{-(2+3\epsilon_1)t} = 2^{(1-3c)t}.
\end{equation}
We also have that $1-1/s < -\epsilon_1 = 1-c$.

Now, for any fixed $(a_0, \mbf a) \in A^t$ satisfying $\Delta_2( F, \la, \psi)\cap J^{(t)}\neq \emptyset$ and any $x\in J^{(t)}$  Taylor's Theorem tells us that there exist points $y_1$ and $y_2$ both lying between $x$ and~$y_0$ satisfying
$$
\F(x) =  \F(y_0) + (x-y_0)\F'(y_0) + (x-y_0)^2 \F''(y_1)/2
$$
and
$$
\F^\prime(x) = \F^\prime(y_0) + (x-y_0)\F^{\prime\prime}(y_2)
$$
respectively. Combining this with inequalities \eqref{bounds}, \eqref{Jsize} and  \eqref{observation}  yields
\begin{eqnarray}
|\F(x)|  & \ll & \psi(2^t) + 2^{-ct}\cdot 2^{(1-1/s)t} +  (2^{-ct})^2\cdot |\mbf a|/2 \nonumber \\ & \ll & 2^{(1-3c) t}+2^{(1-2c)t} +2^{(1-2c)t-1} \nonumber \\& \ll &   2^{(1-2c)t},\label{badinequality1} \
\end{eqnarray}
and also
\begin{equation}\label{badinequality2}
|\F^\prime(x)|\: \: \ll \: \: 2^{(1-1/s)t} + 2^{-ct}\cdot |\mbf a|\: \: \ll \: \: 2^{(1-c) t}+2^{(1-c)t} \: \: \ll \: \: 2^{(1-c)t},
\end{equation}
for any such $(a_0, \mbf a)$ and every $x\in J^{(t)}$.
To conclude our proof  we use similar arguments to those presented in  \cite[pp. 346-351]{DB_inh_dim}. 

\subsubsection{Further auxiliary lemmas and preliminaries}

We begin by listing a collection of useful geometric lemmas from~\cite{DB_inh_dim} that remain unchanged in our setup.
\begin{lem}[Lemma 7, \cite{DB_inh_dim}]\label{DBlem7}
 Consider the plane in $\RR^3$ defined by the equation $Ax_1 + Bx_2 + Cx_3 = D$, where $A,B,C,D$ are integers with $\gcd(A,B,C) = 1$. Then the area $S$ of any
triangle on this plane with integer vertices is at least
$\frac 1 2 \sqrt{A^2 + B^2 + C^2}$.
\end{lem}
Next,  we appeal to a further lemma that follows immediately from arguments presented in~\cite{DB_inh_dim}.

\begin{lem}[\cite{DB_inh_dim}, pp. 346]\label{parallelopiped}
    Fix positive real numbers $t_1, t_2$ and $t_3$ and let $P$ be some affine plane in $\RR^3$ given by the equation $Ax_1+Bx_2+Cx_3=D$. For any functions~$f, \lambda \in C^{(2)}(I)$ and any fixed $x \in I$ consider the parallelepiped $R$ formed by the system of inequalities
    \begin{equation}\label{array1}
    \begin {cases}
    |a_2f(x)+a_1x+a_0+\lambda(x)| \: & \leq \: t_1, \\
    |a_2f^\prime(x)+a_1+\lambda^\prime(x)| \:&\leq\: t_2,\\
    |a_2|\: &\leq\: t_3,
    \end {cases}
    \end{equation}
    where $a_0, a_1, a_2$ are viewed as real variables. Choose a pair of integers $(j, k) \in \left\{1,2,3\right\}$ satisfying $j < k$. Then,   the area $S_{jk}$ of intersection of $P$ with figure defined by the $j$th and $k$th equations in the system~\eqref{array1} satisfies
    $$
    S_{jk} \: \asymp \: \frac{t_j t_k\sqrt{A^2 + B^2 + C^2}}{|T_{jk}|},
    $$
    where
    $$
    T_{12}:=f^\prime(x)(B-Ax)-(C-Af(x)), \quad T_{1,3}:=B-Ax \quad \text{ and } \quad T_{2,3}:=A.
    $$
\end{lem}
    Note that the area of intersection of the plane $P$ with the parallelepiped $R$ as defined in the lemma is not greater than the smallest of~$S_{jk}$.

 Finally, we present a lemma specific to the setting of this paper.
\begin{lem}\label{DBlem9}
    For every $J^{(t)}$ as above, all integer points  $(a_0, \mbf a) \in A^t$
    such that $\Delta_2( F, \la, \psi)\cap J^{(t)}\neq \emptyset$ lie on a single
    affine plane.
\end{lem}

This final lemma is extremely similar to~\cite[Lemma 9]{DB_inh_dim} and the proof follows almost immediately upon replacing the set $A_3(a_0, a_1, a_2)$ from~\cite{DB_inh_dim} with our set $\Delta_2( F_{(a_0, \mbf a)}, \la, \psi)$. One must simply observe that the analogous proof in~\cite{DB_inh_dim} only uses the fact that for all $x\in J^{(t)}$ the function $\F$ satisfies $|\mcal F(x)|\le 2^{(1-2c)t}$ and $|\mcal
F'(x)|\le 2^{(1-c)t}$, and that is precisely what we have shown for our setup via \eqref{badinequality1} and \eqref{badinequality2}.  In fact, the lemma actually holds for all integer triples from the set
$$
A^t_\ast:\: = \: \left\{(a_0, \mbf a) \in \ZZ \times I_2: \: 2^t \leq |\mbf a | < 2^{t+1}, |\mbf a|=r \right\}.
$$
Moreover, note that the discussion of the case of Bad intervals case thus far  has also been independent of the assumption $a_0= a^\prime_0$. (This has actually been true for the entirety of subcase IIb2.) For this reason we will drop this assumption and overestimate the sum $\wdsbbb$ by replacing $A^t$ with~$A^t_\ast$. As we shall see, this will not affect our ability to demonstrate the convergence of $\wdsbbb$.

Upon applying Lemma \ref{DBlem9} it is apparent that there are two final possibilities for each fixed bad interval~$J^{(t)}$:
\begin{enumerate}[i]
\item [(i)] The integer points  $(a_0, \mbf a) \in A_\ast^t$ for which $\Delta_2( F, \la, \psi)\cap J^{(t)}\neq \emptyset$ do not lie on a line.
\item [(ii)] The integer points  $(a_0, \mbf a) \in A_\ast^t$ for which $\Delta_2( F, \la, \psi)\cap J^{(t)}\neq \emptyset$ do lie on a line.
\end{enumerate}

We say  $J^{(t)}$ is a \textit{Bad$(i)$} interval if it satisfies property $(i)$ and a \textit{Bad$(ii)$} interval if it satisfies property $(ii)$.
For one final time we split our volume sum into two cases. Let
\begin{equation*}\label{case_iv3i}
\wdsbbbi=  \displaystyle\sum_{t\, \geq \, \left\lfloor\log_2N \right\rfloor}  \: \: \sum_{\text{Bad$(i)$ intervals } J^{(t)}}\: \:\sum_{(a_0, \mbf a) \: \in \:A^t_\ast}
\: \: |\Delta_2(F,\lambda, \psi) \cap J^{(t)}|^s,
\end{equation*}
and
\begin{equation*}\label{case_iv3ii}
\wdsbbbii=  \displaystyle\sum_{t\, \geq \, \left\lfloor\log_2N \right\rfloor}  \: \: \sum_{\text{Bad$(ii)$ intervals } J^{(t)}}\: \:\sum_{(a_0, \mbf a) \: \in \:A^t_\ast}
\: \: |\Delta_2(F,\lambda, \psi) \cap J^{(t)}|^s.
\end{equation*}
Clearly we have
\begin{equation*}\label{thelastone}
    \wdsbbb \: \ll \: \wdsbbbi \: + \: \wdsbbbii,
\end{equation*}
and so to complete the proof of Theorem \ref{thm6} we must demonstrate that the latter two series converge for $N$ large enough.

\subsubsection{Dealing with Bad$(i)$ intervals}
Let the equation of the plane $P$ in $\RR^3$, on which the points $(a_0, \mbf a) \in A^t$ satisfying $\Delta_2( F, \la, \psi)\cap J^{(t)}\neq \emptyset$ lie, be given by
$Ax_1+Bx_2+Cx_3=D$ for some integers $(A, B, C, D) \in \ZZ^4$.
Recall that for any  $(a_0, \mbf a) \in A^t_\ast$ for which $\Delta_2( F, \la, \psi)\cap J^{(t)}\neq \emptyset$ and any $x \in J^{(t)}$ inequalities \eqref{badinequality1}, \eqref{badinequality2} and
\begin{equation}\label{upper}
|a_2|  \: \asymp \: |\mbf a|\: \ll 2^t
\end{equation}
all hold. Since the integer points $(a_0, \mbf a)$ do not all lie on a line, their maximum number within the parallelepiped $R$ defined by equations \eqref{badinequality1}, \eqref{badinequality2} and \eqref{upper} is bounded above by the maximum number of triangles with integer vertices that may lie in $P \cap R$. It follows upon combining Lemmas~\ref{DBlem7} \& \ref{parallelopiped} that for any Bad$(i)$ interval $J^{(t)}$ this number $N$ of triples $(a_0, \mbf a) \in A^t_\ast$ for which $\Delta_2( F, \la, \psi)\cap J^{(t)}\neq \emptyset$ satisfies
\begin{equation}\label{maxnumber}
N \: \ll \frac{\min\left( S_{12}, \,  S_{13}, \,  S_{23}\right)}{\sqrt{A^2 + B^2 + C^2}} \: \ll \:
\min\left( \frac{2^{(2-3c)t}}{|T_{12}|}, \, \frac{2^{(2-2c)t}}{|T_{13}|}, \, \frac{2^{(2-c)t}}{|T_{23}|}\right).
\end{equation}
However, since $J^{(t)}$ is assumed to be a Bad interval (and so $N > 2^{(3/2 - c - \epsilon_2)t}$) we therefore have for $x \in J^{(t)}$ that
\begin{equation}\label{Tbound1}
    |T_{12}| \, = \, |f^\prime(x)(B-Ax)-(C-Af(x))| \, \ll \, 2^{(1/2-2c+\epsilon_2)t},
\end{equation}
\begin{equation}\label{Tbound2}
|T_{13}| \, = \, |B-Ax| \, \ll \, 2^{(1/2-c+\epsilon_2)t},
\end{equation}
and
\begin{equation}\label{Tbound3}
|T_{23}| \, = \, |A| \, \ll \, 2^{(1/2+\epsilon_2)t}.
\end{equation}
Moreover, combining \eqref{Tbound1} with \eqref{Tbound2} yields
\begin{equation}\label{Tbound4}
|C-Af(x))| \: \leq \: |T_{12}|  + |f^\prime(x)T_{13}|\: \ll \: 2^{(1/2-2c+\epsilon_2)t}+2^{(1/2-c+\epsilon_2)t} \: \ll \: 2^{(1/2-c+\epsilon_2)t}.
\end{equation}

Next, let $\mathcal P_t$ denote the set of all integer triples $(A, B, C)$ that appear as coefficients of a plane associated to some Bad$(i)$ interval~$J^{(t)}$ and by $M_t$ their number; i.e., let $M_t=\# \mathcal P_t$. Then, since any such triple must satisfy \eqref{Tbound2} and \eqref{Tbound4} it follows from the analogous arguments in \cite{DB_inh_dim} that $M_t \: \ll \: |A|$ so long as $\epsilon_2$ is taken to be sufficiently small so that $\epsilon_2 < c-1/2$. In particular, in view of \eqref{Tbound3} and a well known upper bound for the partial sums of the harmonic series we have
\begin{equation}\label{Mbound}
\sum_{(A, B, C) \,\in\, \mathcal P_t}\frac1{|A|^2} \: \:\ll \:\: \sum_{A=1}^{ \left\lfloor 2^{(1/2+\epsilon_2)t}\right\rfloor}\frac1{A}\: \: \:\ll \: \left(\frac 1 2 + \epsilon_2\right)\cdot t \:\: \ll \:\: t.
\end{equation}

Furthermore, inequality \eqref{Tbound2} implies that for any two Bad$(i)$ intervals $J_1^{(t)}$ and $J_2^{(t)}$ whose associated planes have the same coefficients $(A, B, C)$ that
$$
|A(x-y)| \: = \: |(B-Ax)-(B-Ay)|\: \ll \:  2^{(1-2c)t/2},
$$
for  $x\in J_1^{(t)}$ and $y\in J_2^{(t)}$. In turn,
$$
|x-y| \: \ll \: \frac{2^{(1/2-c+\epsilon_2)t}}{|A|}, \qquad x\in J_1^{(t)}, \: \: y\in J_2^{(t)}.
$$
We conclude that for fixed $(A, B, C) \in \mathcal P_t$ a given point~$x$ lying in any set $\Delta_2( F, \la, \psi)$ under consideration may fall only within some interval of length $\ll 2^{(1/2-c+\epsilon_2)t}\,|A|^{-1}$.
Therefore, in view of \eqref{Jsize} the number of Bad$(i)$ intervals associated with the triple
$(A, B, C)\in \mathcal P_t$ is at most  $\ll 2^{(1/2+\epsilon_2)t} \, |A|^{-1}$.


Upon utilising the bounds \eqref{maxnumber}, \eqref{dflam1} and \eqref{Mbound} respectively, we deduce that for $N$ sufficiently large the quantity $\wdsbbbi$ is bounded above by
\begin{eqnarray*}
 && \sum_{t\, \geq \, \left\lfloor\log_2N \right\rfloor}  \: \: \sum_{\text{Bad$(i)$ intervals } J^{(t)} }\: \:\sum_{(a_0, \mbf a) \: \in \:A^t_\ast}
\: \: |\Delta_2(F,\lambda, \psi) \cap J^{(t)}|^s
\\ &\leq&
\sum_{t\, \geq \, \left\lfloor\log_2N \right\rfloor}  \: \:
\frac{\psi(2^t)^{s/2}}{2^{st/2}} \: \:
    \sum_{(A, B, C) \,\in\, \mathcal P_t} \:\:
\sum_{\substack{\text{Bad$(i)$ intervals } J^{(t)} \\ \text{associated with } (A, B, C) }}
\: \:\sum_{\substack{(a_0, \mbf a) \: \in \:A^t_\ast: \\ \Delta_2(F,\lambda, \psi) \cap J^{(t)}\neq \emptyset }} \: \: 1
\\ &\ll&
 \sum_{t\, \geq \, \left\lfloor\log_2N \right\rfloor}  \: \:
 \frac{\psi(2^t)^{s/2}}{2^{st/2}} \: \:
 \sum_{(A, B, C) \,\in\, \mathcal P_t} \:\:
 2^{(1/2+\epsilon_2)t} |A|^{-1} \cdot 2^{(2-c)t}|A|^{-1}
\\ &=&
\sum_{t\, \geq \, \left\lfloor\log_2N \right\rfloor}  \: \: \psi(2^t)^{\frac s 2} \cdot 2^{(\frac32-\frac s2+\epsilon_2-\epsilon_1)t}\sum_{(A, B, C) \,\in\, \mathcal P_t}\: |A|^{-2}
\\ &\ll&
\sum_{t\, \geq \, \left\lfloor\log_2N \right\rfloor}  \: \: \psi(2^t)^{\frac s 2} \cdot 2^{(\frac32-\frac s2+\epsilon_2-\epsilon_1)t} \cdot t
\\&\ll&
 \sum_{t\, \geq \, \left\lfloor\log_2N \right\rfloor}  \: \: \psi(2^t)^{\frac s 2} \cdot 2^{(\frac32-\frac s2+2\epsilon_2-\epsilon_1)t}. \
\end{eqnarray*}
Finally, take $\epsilon_2$ small enough in terms of $\epsilon_1$ so that $\epsilon:=2 \epsilon_1-4\epsilon_2 >0$. Then, by \eqref{psi_low} we have
$$\psi(2^t) \: \geq \: 2^{(1-\frac{3+\epsilon}{s})t},
$$
and in turn
$$
\psi(2^t)^{\frac s 2} \: \:\geq \: \:2^{(\frac s 2-\frac{3+\epsilon}{2})t} \: \:  = \:  \: 2^{(\frac s 2-\frac 3 2 +2\epsilon_2-\epsilon_1)t}.
$$
So, by \eqref{condensation} we deduce
$$
\wdsbbbi \: \: \ll \: \:  \sum_{t\, \geq \, \left\lfloor\log_2N \right\rfloor}  \: \: \psi(2^t)^{s} \: 2^{(3-s)t} \:\:< \: \:\infty,$$
as required.


\subsubsection{Dealing with Bad$(ii)$ intervals}
In this final subcase we consider Bad  intervals~$J^{(t)}$ for which  all integer points  $(a_0, \mbf a) \in A_\ast^t$  satisfying $\Delta_2( F, \la, \psi)\cap J^{(t)}\neq \emptyset$ lie on some affine line in $\RR^3$. For a fixed Bad$(ii)$ interval $J^{(t)}$ let $L= \boldsymbol \alpha +\ell \boldsymbol\beta$ denote this associated line, where $\ell$ is the real parameter and $\boldsymbol\alpha, \boldsymbol\beta \in \RR^3$. Without loss of generality  we may assume that $\boldsymbol\alpha=(\al_0, \al_1, \al_2)$ is an integer vector and $\boldsymbol \beta=(\be_0, \be_1, \be_2)$ is an integer vector connecting $\boldsymbol\alpha$ and the nearest integer point to $\boldsymbol\alpha$  on $L$. Indeed, the line~$L$ must contain integer points by definition. It follows that all vectors $(a_0, \mbf a) \in A_\ast^t$ associated with $J^{(t)}$ must be of the form
$$
(a_0, a_1, a_2) \: = \: (\alpha_0 + k \beta_0, \: \alpha_1 + k \beta_1, \: \alpha_2 + k \beta_2), \qquad k \in \ZZ.
$$
Moreover, since $J^{(t)}$ is a Bad interval there are at least $2^{(3/2 - c - \epsilon_2)t}$ distinct values of $k$ taken. Therefore, there exist $k_1$ and $k_2$ such that $|k_1-k_2|\geq2^{(3/2 - c - \epsilon_2)t}$. Let $(a^\prime_0, a^\prime_1, a^\prime_2) $ and $(a^{\prime\prime}_0, a^{\prime\prime}_1, a^{\prime\prime}_2) $ be the triples associated with $k_1$ and $k_2$ respectively.

Recalling inequalities \eqref{a0counting} and \eqref{upper},  we have for every triple $(a_0, a_1, a_2)$ under consideration that $|a_j| \ll 2^t$ for $j=0,1,2$. It follows that
$$
|\beta_j(k_1-k_2)| \: \: =  \: \: |a^{\prime}_j-a^{\prime\prime}_j|  \:\: \ll \: 2^t, \qquad j=0,1,2,
$$
and in turn that
$$
|\beta_j| \: \: \ll \: \: 2^{(c+\epsilon_2-1/2)t}, \qquad j=0,1,2.
$$
In particular, we have
\begin{equation}\label{betabounds}
\max \left(|\beta_0|, \,  |\beta_1|, \, |\beta_2|\right)   \: \: \ll \: \: 2^{(c+\epsilon_2-1/2)t}.
\end{equation}
Furthermore, inequality \eqref{badinequality1} yields
\begin{equation*}
|(\beta_2 f(x)+\beta_1 x +\beta_0  )(k_1-k_2)| \: \: = \: \:
|\F_{(a^\prime_0, a^\prime_1, a^\prime_2)}(x)-\F_{(a^{\prime\prime}_0, a^{\prime\prime}_1, a^{\prime\prime}_2)}(x)|\:  \: \ll \: \: 2^{(1-2c)t},
\end{equation*}
and inequality \eqref{badinequality2} yields
\begin{equation*}
|( \beta_2 f^\prime(x) +\beta_1 )(k_1-k_2)| \: \: = \: \:
|\F^\prime_{(a^\prime_0, a^\prime_1, a^\prime_2)}(x)-\F^\prime_{(a^{\prime\prime}_0, a^{\prime\prime}_1, a^{\prime\prime}_2)}(x)|\:  \: \ll \: \: 2^{(1-c)t}.
\end{equation*}
In turn, since $|k_1-k_2|\geq2^{(3/2 - c - \epsilon_2)t}$ we have
\begin{equation}\label{betabounds1}
|\beta_2 f(x)+\beta_1 x +\beta_0  | \: \: \ll \: \: 2^{(\epsilon_2-c-1/2)t}  \: \:  \ll \: \: 2^{(\epsilon_2-1/2)t},
\end{equation}
and
\begin{equation}\label{betabounds2}
 | \beta_2 f^\prime(x) +\beta_1 | \: \:  \ll \: \: 2^{(\epsilon_2-1/2)t}.
\end{equation}
 We now appeal to one final lemma.
\begin{lem}[\cite{DB_inh_dim}, pp. 349]\label{cheat}
    Fix any non-degenerate $f \in C^{(2)}$, any $\delta \in (0,1)$ and a positive real parameter~$r$. Then for any fixed $\beta_2 \in \ZZ\setminus \left\{0 \right\}$  the number of integer solutions $(\beta_0, \beta_1, \beta_2)$ to the system of inequalities
    \begin{equation}\label{array}
    \begin {cases}
    |\beta_2 f(x)+\beta_1 x +\beta_0 | \:& \leq \: r^{-\delta}, \\
    |\beta_2 f^\prime(x) +\beta_1| \:& \leq\: r^{-\delta},\\
    \max \left(|\beta_0|, \,  |\beta_1|, \, |\beta_2|\right)\: &  \leq\: r,
    \end {cases}
    \end{equation}
    is $\ll |\beta_2|$.
\end{lem}
The proof of this lemma follows upon minor modification to the argument on page~349 of \cite{DB_inh_dim} and appeals to the change of variables exhibited in \cite[Lemma 6]{DB_inh_dim}. We apply Lemma~\ref{cheat} taking  $r=2^{(c+\epsilon_2-1/2)t}$ and
$$\delta=\frac{1-2\epsilon_2}{2c+2\epsilon_2-1}, \qquad \text{ so that} \qquad  r^{-\delta} \: \: =\: \:  2^{-\frac{(1-2\epsilon_2)(c+\epsilon_2-1/2)t}{2c+2\epsilon_2-1}} \: \: =\: \:  2^{(\epsilon_2-1/2)t}.$$
This yields that for every fixed $\beta_2$ the number of integer points $(\beta_0, \beta_1, \beta_2)$ which appear as parameters in a line $L$ associated to some Bad$(ii)$ interval $J^{(t)}$ is $\ll |\beta_2|$.


Equation~\eqref{betabounds1} holds for all $x$ inside a given 
interval $J^{(t)}$. However, one can find a point $x\in J^{(t)}$ for
which the estimate can be improved. By definition, the Bad interval $J^{(t)}$ must intersect more than $2^{3/2-c-\epsilon_2}$ sets of the form $\Delta_2(F_{(a_0, \mbf a)}, \lambda, \psi)$ for some $(a_0, \mbf a) \in A_\ast^t$. If follows from the pigeon-hole principle that there must exist two numbers
$x_1\in \Delta_2(F_{(a_0', \mbf a')},\lambda, \psi)$ and $x_2\in
\Delta_2(F_{(a''_0, \mbf a'')},\lambda, \psi)$ both satisfying (\ref{bounds}) such that
\begin{equation}\label{thefix}
|x_1-x_2|\le 2^{(-c-3/2+c+\epsilon_2)t} = 2^{(-3/2+\epsilon_2)t}
\end{equation}
for some triples $(a'_0, \mbf a')$ and $(a''_0, \mbf a'')$ from $A_\ast^t$ satisfying $\Delta_2(F_{(a_0', \mbf a')},\lambda, \psi) \cap J^{(t)} \neq \emptyset$ and $\Delta_2(F_{(a_0'', \mbf a'')},\lambda, \psi) \cap J^{(t)}\neq \emptyset$ respectively. Then, by utilising Taylor's formula and all of the inequalities (\ref{bounds}), (\ref{psibound}), (\ref{observation}), (\ref{vb6}) and (\ref{thefix}), one has
$$
|\F_{(a_0',\mbf a')} (x_2)|\ll \psi(2^t) + |x_1-x_2|\cdot 2^{(1-1/s)t}
+\frac{|x_1-x_2|^2}{2} \cdot 2^t \ll 2^{-dt},
$$
where
$$
d:=\min\left\{\frac12 +\frac1s-\epsilon_2, \: 2-2\epsilon_2, \:
\frac3s-1\right\}.
$$
Since $s\in (0,1)$, one can then take $\epsilon_2$ and $\epsilon_1$ so small
that
\begin{equation}\label{eq_dbound}
d\ge 3/2 + 4\epsilon_1+4\epsilon_2.
\end{equation}
By subtracting $\F_{(a_0',\mbf a')}(x_2)$ from $\F_{(a''_0,\mbf
	a'')}(x_2)$ we get from (\ref{bounds}) and (\ref{psibound}) that
\begin{eqnarray}
|\beta_2f(x_2)+\beta_1x_2+\beta_0| \: \ll \: |\F_{(a^\prime_0, \mbf a')}(x_2)-\F_{(a^{\prime\prime}_0, \mbf a'')}(x_2)| & \ll& 2^{-dt} + \psi(2^t) \nonumber \\
& \ll & 2^{-dt} + 2^{(1-3c)t} \nonumber \\
& \ll & 2^{-dt}. \label{ineq_fline} \
\end{eqnarray}
In other words, for each Bad$(ii)$ interval $J^{(t)}$ there must
exist at least one number $x$ inside it which
satisfies~\eqref{ineq_fline}.
Now, consider two Bad$(ii)$ intervals~$J_1^{(t)}$ and $J_2^{(t)}$
such that their associated lines $L_1$ and $L_2$ share the same
parameters $(\beta_0, \beta_1, \beta_2)$. Then there are
numbers $x \in J_1^{(t)} $ and $y \in J_2^{(t)}$ such that they both
satisfy estimate~\eqref{ineq_fline}. Since
$|\F_{\beta_0,\beta_1,\beta_2}''(x)| = |\beta_2f''(x)|\gg
|\beta_2|$, by Lemma \ref{Pya_lem} the measure of numbers satisfying it is
at most
$$
(2^{dt}|\beta_2|)^{-1/2}.
$$
Since $\F_{\beta_0,\beta_1,\beta_2}''(x)$ does not change sign,
the set of numbers $x$ satisfying~\eqref{ineq_fline} is the union of
at most two intervals and therefore the number of Bad$(ii)$
intervals $J^{(t)}$ whose associated line has parameters $(\beta_0,
\beta_1, \beta_2)$ is \begin{equation}
\label{anothercount}\ll 2^{(c-d/2)t}|\beta_2|^{-1/2}.
\end{equation}

Finally, let $\mathcal{L}_t$ denote the set of all integer  triples
$(\beta_0, \beta_1, \beta_2)$ which appear as parameters in a line
$L$ associated to some Bad$(ii)$ interval $J^{(t)}$. Observe that
$|\alpha_2+k\beta_2| \ll 2^t$ for every $k$ associated to one of the
triples~$(a_0, \mbf a)$ under consideration. It follows  that for a
fixed integer point $(\beta_0, \beta_1, \beta_2) \in \mathcal{L}_t$
there are at most $\ll 2^t|\beta_2|^{-1}$ triples $(a_0, \mbf a) \:
\in \:A^t_\ast$ for which $\Delta_2(F,\lambda, \psi) \cap J^{(t)}\neq
\emptyset$. Also, upon applying Lemma~\ref{cheat} and utilising the
bound~\eqref{betabounds} we have
\begin{equation}\label{Mbound1}
\sum_{(\beta_0, \beta_1, \beta_2) \,\in\, \mathcal
L_t}\frac1{|\beta_2|^{3/2}} \: \:\ll \:\: \sum_{\beta_2=1}^{
\left\lfloor
2^{(c+\epsilon_2-1/2)t}\right\rfloor}\frac1{|\beta_2|^{1/2}}\: \:
\:\ll \: 2^{(c/2+\epsilon_2/2-1/4)t}.
\end{equation}

In view of the above discussion, it follows from
\eqref{dflam1}, \eqref{anothercount},~\eqref{Mbound1}, ~\eqref{psi_low} and (\ref{eq_dbound}) respectively
that for a large enough~$N$ and a small enough $\epsilon_1$ and
$\epsilon_2$ we have
\begin{eqnarray*}
    \wdsbbbii &=& \sum_{t\, \geq \, \left\lfloor\log_2N \right\rfloor}  \: \: \sum_{\text{Bad$(ii)$ intervals } J^{(t)} }\: \:\sum_{(a_0, \mbf a) \: \in \:A^t_\ast}
    \: \: |\Delta_2(F,\lambda, \psi) \cap J^{(t)}|^s
    \\ &\ll&
    \sum_{t\, \geq \, \left\lfloor\log_2N \right\rfloor}  \: \:
    \frac{\psi(2^t)^{s/2}}{2^{st/2}} \: \:
    \sum_{(\beta_0, \beta_1, \beta_2) \,\in\, \mathcal L_t} \:\:
        \frac{2^t}{|\beta_2|} \: \: \sum_{\substack{\text{Bad$(ii)$ intervals } J^{(t)} \\ \text{associated with } (\beta_0, \beta_1, \beta_2) }}
    \: \: 1
    \\ &\ll&
    \sum_{t\, \geq \, \left\lfloor\log_2N \right\rfloor}  \: \:
    \frac{\psi(2^t)^{s/2}}{2^{st/2}} \: \:
    \sum_{(\beta_0, \beta_1, \beta_2) \,\in\, \mathcal L_t} \:\:
        2^t |\beta_2|^{-1} \cdot 2^{(c-d/2)t}|\beta_2|^{-1/2}
    \\ &=&
    \sum_{t\, \geq \, \left\lfloor\log_2N \right\rfloor}  \: \: \psi(2^t)^{\frac s 2} \cdot 2^{t(\frac54-\frac s2-\epsilon_1-2\epsilon_2)}\sum_{(\beta_0, \beta_1, \beta_2) \,\in\, \mathcal L_t}\: |\beta_2|^{-\frac 3 2}
    \\ &\ll&
    \sum_{t\, \geq \, \left\lfloor\log_2N \right\rfloor}  \: \: \psi(2^t)^{\frac s 2} \cdot 2^{(\frac32-\frac s2-\frac 1 2 \epsilon_1-\frac 3 2\epsilon_2)t}
    \\&\ll&
    \sum_{t\, \geq \, \left\lfloor\log_2N \right\rfloor}  \: \:
    2^{-(\frac{1}{2}\epsilon_1+\frac32\epsilon_2)t}
    \qquad < \quad \infty,\
%
%
%
%
\end{eqnarray*}
and the proof is complete. 

\

{\footnotesize

\bibliographystyle{abbrv}
\bibliography{Mumtaz}}

\end{document}